\documentclass[11pt]{article}
\usepackage{latexsym,epsfig}
\font\smf=cmr7
\def\bx{\hfill $\Box$}
\def\cc{{\cal C}}
\def\ll{{\cal L}}
\def\prc{\mathrm{Prob}(\cc)}
\def\prl{\mathrm{Prob}(\ll)} 
\def\eqdef{\, =\kern -12pt\raise 6pt\hbox{{\smf def}}\, }
\newtheorem{theorem}{Theorem}
\newtheorem{lemma}{Lemma}
\newtheorem{proposition}{Proposition}

\newtheorem{corollary}{Corollary}

\font\smf=cmr7

\begin{document}
{\topskip 1in
\begin{center}
{\huge
The distributions of the entries of Young tableaux}
\vspace{50pt}

{\Large
Brendan D. McKay\kern1pt\footnote{Dept. of Computer Science, Australian
National University, 
ACT 0200, Australia; 
\texttt{\small <bdm@cs.anu.edu.au>}}, 
Jennifer Morse\kern1pt\footnote{Dept. of Mathematics, University of
Pennsylvania, 
Philadelphia, PA 19104-6395;
\texttt{\small <morsej@math.upenn.edu>}
}, and Herbert S. Wilf\kern2pt\footnote{Dept. of Mathematics, University of
Pennsylvania, 
Philadelphia, PA 19104-6395;
\texttt{\small <wilf@math.upenn.edu>}
}}
\end{center}

\vspace{.8in}

\begin{abstract}
Let $T$ be a standard Young tableau of shape
$\lambda\vdash k$. 
We show that the probability that a randomly chosen
Young tableau of $n$ cells contains $T$ as a subtableau is, 
in the limit $n\to\infty$, equal to $f^{\lambda}/k!$, 
where $f^\lambda$ is the number of all tableaux of 
shape $\lambda$. In other words, 
\textit{the probability that a large tableau contains $T$ 
is equal to the number of tableaux whose shape is that of $T$, 
divided by~$k!$}.

We give several applications, 
to the probabilities that a set of prescribed entries will appear in a 
set of prescribed cells of a tableau, and to the probabilities that 
subtableaux of given shapes will occur. 

Our argument rests on a notion of quasirandomness of families of permutations,
and we give sufficient conditions for this to hold.
\end{abstract}
}

\vspace{.6in}

\newpage
\section{Main results}
\label{sec:result}
Our basic result is the following.
\begin{theorem}
\label{th:basic}
Fix a standard Young tableau $T$ of shape $\lambda\vdash k$, 
let $N(n;T)$ be the number of tableaux of $n$ cells that contain $T$ 
as a subtableau,\footnote{A subtableau of a tableau $T$ of $n$ cells is a
tableau that is formed by the letters $1,2,\dots,k$ in $T$, for some $k\le
n$.}
and let $t_n$ be the number of all tableaux of $n$ cells. 
Then we have
\begin{equation}
\label{eq:basic}
\lim_{n\to\infty}\frac{N(n;T)}{t_n}=\frac{f^{\lambda}}{k!},
\end{equation}
where $f^\lambda$ is the number of all tableaux of shape $\lambda$. 
In other words, the probability that a large tableau contains $T$ 
is equal to the number of tableaux whose shape is that of $T$, 
divided by $k!$.
\end{theorem}

We now state two corollaries of this theorem, 
after which we will discuss several applications.
Two excellent references regarding the general theory of 
tableaux are \cite{fu} and \cite{kn}.

\begin{corollary}
\label{cor:first}
Let $\cc$ be a collection of Young tableaux, none of which is a 
subtableau of any other in the collection, and
let $N(n;\cc)$ be the number of Young tableaux of $n$ cells 
which have a subtableau in $\cc$.
The probability that a randomly chosen tableau of $n$ cells has
a subtableau in $\cc$ is then $N(n;\cc)/t_n$, 
where $t_n$ is the number of tableaux of $n$ cells 
(equivalently, the number of involutions of $n$ letters). 
We have
\begin{equation}
\label{eq:prdef}
\prc\eqdef 
\lim_{n\to\infty}\frac{N(n;\cc)}{t_n}=
\sum_{T\in\,\cc}\frac{f^{\lambda(T)}}{|T|!}\,,
\end{equation}
where $\lambda(T)$ is the shape of tableau $T$
and $|T|$ is the number of cells in $T$.
\end{corollary}
\noindent Thus we can speak of ``the probability that a Young tableau 
has a subtableau appearing in $\cc$,"  without reference to the size, $n$, 
of the tableau.  This phrase will mean the limit in (\ref{eq:prdef}).

\smallskip

The next corollary is the special case of 
Corollary \ref{cor:first} in which the distinguished 
list $\cc$ of tableaux is defined by a list of allowable shapes. 
\begin{corollary}
\label{cor:second}
Let $\ll$ be a list of Ferrers diagrams with no shape
a subshape of another in the list,
and let $N(n;\ll)$ be the number of Young tableaux of 
$n$ cells which have a subtableau with shape in $\ll$.
The probability that a tableau of $n$ cells has 
such a subtableau is then $N(n;\ll)/t_n$, and we have
\begin{equation}
\label{eq:prldef}
\prl\eqdef \lim_{n\to\infty}\frac{N(n;\ll)}{t_n}=
\sum_{\lambda\in\ll}\frac{(f^{\lambda})^2}{|\lambda|!}.
\end{equation}
\end{corollary}

{From} these results, we will deduce a number of interesting consequences:
\begin{enumerate}
\item Let $\cc$ be the list of all tableaux of $k$ cells such that the 
letter $k$ lives in the $(i,j)$ position, for some fixed $(i,j)$. Then $\prc$ 
is the probability that a Young tableau has the entry $k$ in its $(i,j)$ 
position. We will find a rather explicit formula (see subsection 
\ref{subsec:onecell} below) for this probability. This formula 
was previously found by Regev \cite{re}.
\item Let $\cc$ be the list of all tableaux of $k$ cells in which a 
certain fixed collection of cells contain prescribed entries. Then $\prc$ 
is the probability that a Young tableau has the prescribed entries in 
the prescribed cells. We will find a rather explicit formula for this 
probability (see subsection \ref{subsec:several} below). In 
the case where the fixed collection of cells consists of just two cells, this 
formula was also previously found by Regev \cite{re}, who also found this 
two-cell probability with a variety of measures on the space of tableaux. 
Our result, while it applies to arbitrary collections of prescribed cells, 
holds only in the uniform measure on tableaux.
\item Let $\ll$ be the list of all partitions of the integer $k$ whose 
parts are $\le 2$, in Corollary \ref{cor:second}. Then $\prl$ is 
the probability that a Young tableau has its smallest $k$ letters in 
just two columns, and we'll find an explicit formula for it 
(see subsection \ref{subsec:twocols} below).
\end{enumerate}

Finally, in section \ref{sec:cell12}, we will find the probability that 
the $(1,2)$ entry of a tableau of $n$ cells is $k$, in the form of an 
exact formula that is valid for every $n,k$. The asymptotic form of this 
result will illustrate the rate of approach to the limit in the more 
general theorems already cited above.

\bigskip

We thank the anonymous referee who noted that
our proof of Corollaries \ref{cor:first} and \ref{cor:second} actually proved
them in the form shown 
here, which is more general than our original statement.

\section{Proof of Theorem \ref{th:basic}}
The set of
letters $\{1,2,\ldots,\ell\}$ is denoted $[\ell]$. We begin with a small
observation. 
\begin{proposition}
In the Robinson-Schensted (RS) correspondence between 
involutions $\phi$, of $n$ letters, and tableaux $T$, of $n$ cells, 
the subtableau of $T$ in the letters $[k]$ depends 
only on the order of the first $k$ letters in the involution $\phi$, 
and does not depend on their preimages or on the disposition 
of the remaining $n-k$ letters.
\end{proposition}
To see this, note that when a letter $>k$ is inserted into some stage of the 
RS algorithm it cannot disturb the position of any letter $\le k$. \bx

\bigskip

Fix a tableau $T$, of $k$ cells. 
How many involutions of $n$ letters correspond to a tableau 
that contains $T$ as a subtableau? To answer this, 
let $Z(k)$ denote the set of all permutations of $k$ letters 
which correspond, under the RS correspondence, to an ordered pair of 
tableaux $(T,T')$ for some tableau $T'$.  Then an involution of $n$ 
letters will 
correspond to a tableau that contains $T$ iff the set of letters 
$1,2,\dots,k$ in its value sequence appear in one of the arrangements 
in $Z(k)$.

\smallskip

Thus if $\sigma$ is some permutation of $k$ letters, and $F_n(\sigma)$ 
denotes the number of involutions of $n$ letters which contain $\sigma$ 
as a subsequence, then exactly
\[\sum_{\sigma\in Z(k)}F_n(\sigma)\]
involutions of $n$ letters correspond to tableaux which contain $T$ 
as a subtableau, so the probability that a random $n$-tableau 
contains $T$ is
\begin{equation}
\label{eq:prob1}
\sum_{\sigma\in Z(k)}\frac{F_n(\sigma)}{t_n},
\end{equation}
where $t_n$ is the number of $n$-involutions.
What can be said about the summand $F_n(\sigma)/t_n$? It is the probability 
that a random \textit{involution} of $[n]$ contains the letters $1,2,\dots,k$ 
in some particular order $\sigma$. If instead we had wanted the probability 
that a random \textit{permutation} of $[n]$ contains the 
letters $1,2,\dots,k$ in some particular order, the question would have 
been trivial: the required probability would be exactly $1/k!$, 
no matter what the ``particular order'' was.

\smallskip

We claim that for involutions the answer is essentially the same, up to a
term 
that is $o(1)$ as $n\to\infty$. 
\begin{lemma}
\label{lem:typ}
Let $\sigma$ be a fixed permutation of $k$ letters. 
The probability that a random involution of $n$ letters contains 
$\sigma$ as a subsequence is $1/k!+o(1)$, for $n\to\infty$.
\end{lemma}
We will prove this lemma in the next section as a corollary of a more general 
theorem about the quasirandomness of families of permutations.

\smallskip

However, for the moment let us imagine that we have proved the Lemma, and 
we will now finish the proof of Theorem \ref{th:basic}. 
By (\ref{eq:prob1}) and the Lemma, the probability that a tableau 
of $n$ letters contains a given subtableau $T$ of $k$ letters is
\[\sum_{\sigma\in Z(k)}\left(\frac{1}{k!}+o(1)\right)=\frac{|Z(k)|}{k!}+o(1)
\qquad (n\to\infty).\]
Since $Z(k)$ is the number of all permutations of $k$ letters 
corresponding to ordered pairs of the form $(T,T')$ for some $T'$,
well-known RS theory gives that this is simply
the number of tableaux $T'$ whose shape is that of $T$,
i.e. $f^{\lambda(T)}$.
Thus the probability that a tableau of $n$ letters 
contains a fixed $T$ of $k$ letters as a subtableau is
\[\frac{f^{\lambda(T)}}{k!}+o(1)\qquad (n\to\infty),\]
 and the proof of Theorem \ref{th:basic} is complete. \bx

\bigskip

Corollary \ref{cor:first} follows from the theorem by summing over 
the tableaux in the list $\cc$, since two tableaux cannot be subtableaux
of the same larger tableaux unless one is a subtableau of the other.
Corollary \ref{cor:second} follows 
from Corollary \ref{cor:first} since, if $\cc$ is the 
list of all of the tableaux whose shapes are in~$\ll$,
\begin{eqnarray*}
\sum_{T\in\cc}
\frac{f^{\lambda(T)}}{|T|!}
&=&
\sum_{\lambda\in\ll}\ \sum_{\{T:\lambda(T)=\lambda\}}
\frac{f^{\lambda(T)}}{|T|!} \\
&=&\sum_{\lambda\in\ll}\ 
\sum_{\{T:\lambda(T)=\lambda\}}
\frac{f^{\lambda}}{|\lambda|!} \\
&=&
\sum_{\lambda\in\ll}
\frac{f^{\lambda}}{|\lambda|!} 
\sum_{\{T:\lambda(T)=\lambda\}}\kern-12pt 1
\quad=\quad
\sum_{\lambda\in\ll}\frac{(f^\lambda)^2}{|\lambda|!}
\, .
\end{eqnarray*}
\section{Involutions are typical}
\label{sec:qrnd}

\def\Pn{{\cal P}_n}

In this section we will prove a proposition that implies Lemma \ref{lem:typ}
above.

Let $\cal P$ be a collection of permutations such that
$\Pn={\cal P}\cap {\cal S}_n$ is non-empty for infinitely many values of~$n$,
where ${\cal S}_n$ is the set of all permutations of $[n]$.

If $\tau$ is a sequence of $k$ distinct elements of $[n]$,
let $h(n,\tau)$ be the number of elements of $\Pn$ that have $\tau$ 
as a subsequence.  If $\Pn\ne\emptyset$, the
probability that a random element of $\Pn$ has $\tau$ as a
subsequence is $\tilde p(n,\tau) = h(n,\tau)/|\Pn|$.

Inspired by the terminology of Chung and Graham \cite{cgw},
we say that $\cal P$ is {\it quasirandom\/} if, for each $k\ge 1$,
\[\lim_{n\to\infty}\max_\tau
  \left|\, \tilde p(n,\tau)-{1\over k!} \right| \to 0,\]
where the limit is restricted to those $n$ for which $\Pn$ is
nonempty and the maximum is over all sequences $\tau$ of $k$ distinct
elements of $[n]$.

In this section we will first give a general criterion, involving the
fixed points of the permutations in the family ${\cal P}$, that guarantees
the quasirandomness of the family. Then we will show that the involutions
satisfy this criterion, which is the result that we need for the analysis
of the limiting distributions of the entries of standard tableaux.

\begin{theorem}\label{thm:quasi}
If each $\Pn$ is a union of conjugacy classes of ${\cal S}_n$, and the
average number of fixed points of elements of $\Pn$ is $o(n)$,
then $\cal P$ is quasirandom.
\end{theorem}

\noindent{\bf Proof.}\,
We restrict $n$ to values for which $\Pn\ne\emptyset$ and fix
$k\ge 1$.
Let $I$ be any $k$ subset of $[n]$, and let
${\cal S}$ be the set of all permutations of~$I$.
Also let $f_n$ be the average number of fixed points of elements of~$\Pn$.

The set $\Pn$ can be expressed as a disjoint union
\[\Pn = A(I) \cup \; \bigcup_{\tau\in {\cal S}} B(\tau),\]
where $A(I) =\{\phi\in\Pn \,|\, \phi(I)\cap I\ne\emptyset\}$, and
\[B(\tau) = \{\phi\in\Pn\,|\, \phi\notin A(I) \,
   \hbox{and $\phi$ contains $\tau$ as a subsequence}\}.\]
The basic idea of the proof is that $A(I)$ is small compared to
$\Pn$ and the size of $B(\tau)$ is independent of~$\tau$.

We begin by showing that $A(I)$ is small.  
Since $\Pn$ is closed under conjugation, all elements of $[n]$ are equally 
likely to be fixed points of members of $\Pn$. 
Thus if we let $t\in I$, then the probability that a random element of $\Pn$ 
fixes $t$ is exactly $f_n/n$, For the same reason, the probability 
that $t$ is mapped onto a specified element of $I$ other than $t$ is exactly
\[1-f_n/n\over n-1.\]
Therefore, the probability that
$t$ is mapped to an element of $I$ is
\[{f_n\over n} + (k-1){1-f_n/n\over n-1} = {n(k-1)+f_n(n-k)\over n(n-1)}.\]
Consequently, the probability $q(I)$ that a random element of $\Pn$
is in $A(I)$ is
\begin{equation}
\label{eq:qest}
q(I) = {|A(I)|\over |\Pn|} \le k\Bigl({n(k-1)+f_n(n-k)\over n(n-1)}\Bigr),
\end{equation}
which shows that $q(I)=o(1)$ if $f_n=o(n)$.

Next, let $\tau_1$ and $\tau_2$ be permutations of~$I$
(that is, elements of ${\cal S}$).
Let $\phi$ be the element of ${\cal S}_n$ which fixes everything 
not in $I$ and maps $\tau_1$ element-wise onto $\tau_2$.
Then conjugation by $\phi$ is a bijection from $B(\tau_1)$ to
$B(\tau_2)$, so $B(\tau_1)$ and $B(\tau_2)$ have the same size.
Thus the sets $B(\tau)$ have the same size for any
$\tau\in{\cal S}$.

Altogether, then, we find that
\[{1-q(I)\over k!} \le \tilde p(n,\tau) \le {1-q(I)\over k!}+q(I),\]
which is the result we want.  The left and right sides come from
supposing that none or all of the elements of $A(I)$, respectively,
contain the subsequence $\tau\in{\cal S}$. \bx

To apply Theorem~\ref{thm:quasi} to the set of all involutions, it
suffices to show that involutions have on average $o(n)$ fixed points.
Since the involutions fixing some specified point are just the
involutions of the remaining points, we have that the average
number of fixed points is exactly $nt_{n-1}/t_n$, which, in view of the well
known asymptotic behavior of $t_n$, viz.
\begin{equation}\label{eq:invol}
t_n=\frac{1}{\sqrt{2}}n^{n/2}
\exp{\left(-\frac{n}{2}+\sqrt{n}-\frac{1}{4}\right)}(1+o(1))\qquad 
(n\to\infty),
\end{equation}
 is
$(1+o(1))\sqrt n$. \bx

\section{Applications}
In this section we will apply Theorem \ref{th:basic} to the examples that
were 
listed above in section \ref{sec:result}.
\subsection{Occupancy of a cell in a tableau}
\label{subsec:onecell}
For fixed positive integers $(i,j)$ and $k$, what is the probability that a 
Young tableau has its $(i,j)$ entry equal to $k$?

The list $\cc$ here consists of all tableaux of $k$ cells whose $(i,j)$ entry 
is $k$. The required probability is, by Corollary \ref{cor:first}, 
\begin{equation}
\label{eq:oneentry}
\frac{1}{k!}\sum_{T_{i,j}=k;\, |T|=k}f^{\lambda(T)}=
\frac{1}{k!}{\sum_{|\lambda|=k}^{}}\kern-2pt\raise7pt\hbox{$\prime$}
f^{\lambda}f^{\lambda-(i,j)}\,,
\end{equation}
where the latter sum extends over all partitions $\lambda$ of the integer $k$ 
whose Ferrers diagram has the cell $(i,j)$ as a corner position, 
and $\lambda-(i,j)$ is the Ferrers diagram of $\lambda$ after removing 
the corner $(i,j)$. 

In particular cases one can make this quite explicit. 
Short computations with the hook formula now reveal, 
for instance, the following.
\begin{enumerate}
\item \label{itm:first} The probability that the $(1,2)$ entry of a 
Young tableau is $k$ is $(k-1)/k!$, for $k=2,3,\dots$.
\item The probability that the $(1,3)$ entry of a Young tableau is $k$ is 
\[\frac{(2k-2)!}{(k-3)!\,k!\,(k+1)!}\qquad   (k=3,4,5,\dots),\]
and the probability that the $(2,2)$ entry is $k+1$ is exactly the same! 
These cases were previously derived by Regev \cite{re}, and the fact that 
the $(1,3)$ and the $(2,2)$ answers are so related is explained there in a 
more combinatorial way.
\end{enumerate}
We show below a short table of the limiting probability (\ref{eq:oneentry}) 
that the $(i,j)$ entry of a Young tableau is equal to $k$.
\[
\begin{array}{ccccccccccccc}
\mathbf{k\backslash(i,j):}&\mathbf{(1,2)}&\mathbf{(1,3)}&\mathbf{(1,4)}&
\mathbf{(1,5)}&\mathbf{(1,6)}&\mathbf{(2,2)}&\mathbf{(2,3)}\\
&&&&&&\\
\mathbf{2:}&\frac{1}{2}&0&0&0&0&0&0\\
&&&&&&\\
\mathbf{3:}&\frac{1}{3}&\frac{1}{6}&0&0&0&0&0\\
&&&&&&\\
\mathbf{4:}&\frac{1}{8}&\frac{1}{4}&\frac{1}{24}&0&0&\frac{1}{6}&0\\
&&&&&&\\
\mathbf{5:}&\frac{1}{30}&\frac{7}{30}&\frac{1}{10}&\frac{1}{120}&0&\frac{1}{
4}&0\\
&&&&&&\\
\mathbf{6:}&\frac{1}{144}&\frac{1}{6}&\frac{7}{48}&\frac{1}{36}&\frac{1}{720
}&\frac{7}{30}&\frac{5}{144}&&&&&
\end{array}
\]

\bigskip

Professor Okounkov has kindly communicated to us  (p.c.) another,
independent proof of  the result (\ref{eq:oneentry}).

\subsection{Occupancies of several cells in a tableau}
\label{subsec:several}
Suppose we're given a finite collection of cells $(i_1,j_1),\dots, (i_m,j_m)$ 
and a collection of entries $k_1,\dots,k_m$. Let $K=\max_i{k_i}$. 
In Corollary \ref{cor:first} let the list $\cc$ consist of all tableaux of
$K$ 
cells that have the given entries in the given cells. Then 
the probability that a Young tableau has all of the entries 
\[\{k_r\in\  \mathrm{cell}\ (i_r,j_r): r=1,\dots,m\}\]
 is
\begin{equation}
\label{eq:setofcells}
\frac{1}{K!}\sum f^{\lambda(T)}\,
\end{equation}
where the sum extends over all tableaux of $K$ cells that have the given
set of
entries in the given set of cells.

We show below a short table of the joint distribution of cells $(1,2)$ and
$(1,3)$. That is, the entry in row $r$ and column $s$ below is the limiting
probability (\ref{eq:setofcells}) that a Young tableau $T$ will  have
$T(1,2)=r$ and $T(1,3)=s$.
\[
\begin{array}{ccccccccccccc}
\mathbf{r\backslash s:}
&\mathbf{3}&\mathbf{4}&\mathbf{5}&\mathbf{6}&\mathbf{7}&\mathbf{8}&\mathbf{9}\\
&&&&&&\\
\mathbf{2:}&\frac{1}{6}&\frac{1}{8}&\frac{11}{120}&\frac{7}{120}&\frac{9}
{280}&\frac{209}{13440}&\frac{2431}{362880}\\
&&&&&&\\
\mathbf{3:}&0&\frac{1}{8}&\frac{11}{120}&\frac{7}{120}&\frac{9}{280}&
\frac{209}{13440}&\frac{2431}{362880}\\
&&&&&&\\
\mathbf{4:}&0&0&\frac{1}{20}&\frac{13}{360}&\frac{53}{2520}&\frac{47}{4480}&
\frac{209}{45360}\\
&&&&&&\\
\mathbf{5:}&0&0&0&\frac{1}{72}&\frac{5}{504}&\frac{73}{13440}&\frac{913}{362
880}\\
&&&&&&\\
\mathbf{6:}&0&0&0&0&\frac{1}{336}&\frac{17}{8064}&\frac{5}{4536}&&&&&
\end{array}
\]

\subsection{An interesting special case}
\label{subsec:twocols}
When the distinguished list $\cc$ of $k$-tableaux consists of \textit{all} 
tableaux that have a certain specified list $\ll$ of shapes, or 
partitions of the integer $k$, we can use Corollary \ref{cor:second}. 
It tells us that the probability that in a Young tableau, the subtableau 
formed by the letters $\{1,2,\dots,k\}$ has one of the shapes in a 
given list $\ll$ of shapes is
\[\frac{1}{k!}\sum_{\lambda\in\ll}(f^{\lambda})^2.\]
This will be recognized as a partial sum of the ``Parseval identity''
\begin{equation}
\label{eq:prsvl}
\frac{1}{k!}\sum_{|\lambda|=k}(f^{\lambda})^2=1
\end{equation}
that holds in the symmetric group ${\cal S}_k$. 
In fact Corollary \ref{cor:second} shows that the 
quantity $(f^{\lambda})^2/k!$ is the probability that  in a 
large tableau $T$ the letters $1,2,\dots,k$ will be arranged in the 
shape $\lambda$, and from this interpetation, (\ref{eq:prsvl}) is obvious. 

An example of this type, i.e., where membership in the distinguished list
$\cc$
depends only on the shape of the tableau, is the following: what is the
probability that in a large Young tableau the letters $\{1,2,\dots,k\}$ are
contained in a subtableau of at most two columns?

{From} Corollary \ref{cor:second} and the hook formula it is a brief exercise
to verify that this probability is
\[\frac{1}{(k+1)!}{2k\choose k}\quad (k=1,2,3,\dots).\]

\section{An exact solution for the cell $(1,2)$}
\label{sec:cell12}

One of the consequences of our main theorem is, as we saw in item
\ref{itm:first} of subsection \ref{subsec:onecell} above, the fact that the
probability that the entry in the $(1,2)$ position of a tableau is equal to
$k$
is $(k-1)/k!$. It is instructive to work out this case also independently of
our main theorem because it happens that an exact solution can be found for
each $n$, instead of only a solution for the limiting probability. This will
shed some light on the rate of convergence to the limit (see
(\ref{eq:secondterm}) below).

Let $f(n,k)$ be the number of standard Young tableaux 
$T$ of $n$ cells (and more than one column)  for which $T_{1,2}=k$.  
That is, $f(n,k)$ is the number of standard tableaux with the letter $k$
occuring 
in the first row and second column.
If we look at the first few 
values of $n$ we see the values of $f(n,k)\}_{k=1}^n$ that are shown below.
\begin{eqnarray}
\label{eq:tabl}
n=2:&\quad&0,1\nonumber\\
n=3:&\quad&0,2,1\nonumber\\
n=4:&\quad&0,5,3,1\nonumber\\
n=5:&\quad&0,13,8,3,1\nonumber\\
n=6:&\quad&0,38,24,9,3,1\nonumber\\
n=7:&\quad&0,116,74,28,9,3,1\nonumber\\
n=8:&\quad&0,382,246,93,29,9,3,1\nonumber\\
n=9:&\quad&0,1310,848,321,98,29,9,3,1\nonumber\\
\end{eqnarray}

\smallskip

Based on the Robinson-Schensted correspondence, we will find an 
exact formula (see (\ref{eq:fform}) below) for these numbers. 
\begin{lemma}
For $1\le k\le n-1$, let $F(n,k)$ denote the set of involutions of\/ 
$[n]$ that contain the subsequence $1\,2 \ldots k$ and
in which the letter $k+1$ occurs before $k$. Then
\[f(n,k)=\cases{ 1,&if\/ $k=0$;\cr\noalign{\vskip 3pt}
                |F(n,k-1)|,&if\/ $1\le k\le n$.} \]
\end{lemma}

\noindent{\bf Proof.}\,
Since the RS-correspondence is a bijection between all
involutions and all pairs of standard tableaux $(P,P)$,
it suffices to show that under the RS-correspondence,
an element $\omega \in F(n,k-1)$ is sent to a pair $(T,T)$ where
the first $k-1$ entries in the first column of
$T$ are $1,2,\dots,k-1$ and $k\in T_{1,2}$.
This follows from the definition of column insertion.
That is, a letter $\ell$ is bumped out of the first column only
if a letter $j$, where $j<\ell$, is column inserted after
$\ell$.  If $\ell$ is any of the letters $1,2,\ldots,k-1$,
there is no letter $j$ where $j<\ell$ inserted after
$\ell$ and thus these letters all remain in the first column.
On the other hand, since $k$ is inserted before $k-1$, it 
must be bumped to the second column and takes position $T_{1,2}$. \bx

\bigskip

\begin{lemma}
\label{lem:rcrr}
Let $G(n,k)$ be the number of involutions of 
$[n]$ containing the subsequence $1\,2\ldots k$. Then we have
\begin{equation}
\label{eq:gfrel}
G(n,k) = \sum_{j\ge k} f(n,j+1)\, .
\end{equation}
\end{lemma}

\noindent {\bf Proof.}\,
The set of all involutions containing
the subsequence $1\,2\ldots k$
can be divided into the subset
where $k+1$ occurs before $k$ and the subset
where $k+1$ occurs after $k$.
By definition, this is to say that
$G(n,k)=|F(n,k)|+G(n,k+1)$, which, after summation on $k$, establishes the 
result. \bx

\smallskip

Now we can find an explicit formula for $f(n,k)$ in terms of the 
number of certain involutions.
\begin{theorem}
\label{th:fform}
Let $f(n,k)$ denote the number of standard tableaux on $n$ letters with 
the entry $k$ occuring in the (1,2) position. 
We then have the exact formula
\begin{equation}
\label{eq:fform}
f(n,k)=\sum_{j=0}^{k-1}{n-k\choose k-j-1}t_{n-2k+j+2}-{n-k\choose
k-1}t_{n-2k+1}-{n-k\choose k}t_{n-2k}
\,.
\end{equation}
\end{theorem}

\noindent {\bf Proof.}\,
Suppose $G(n,k)$ is, as in Lemma \ref{lem:rcrr}, the number of involutions 
of $[n]$ which contain the subsequence $1\,2\dots\,k$. The number 
of these whose value at 1 is 1 is $G(n-1,k-1)$. Now consider such an 
involution $\phi$ for which $\phi(1)>1$. Then in fact $\phi(1)>k$. Hence we 
can choose the locations of the subsequence $1\,2\dots\,k$ in ${n-k\choose 
k}$ ways. Having done that, the values of $\phi$ at $1,2,\dots ,k$ are also 
determined since an involution is composed only of elementary transpositions. 
That leaves $n-2k$ values, which can be any involution of 
$n-2k$ letters. Thus
\[G(n,k)=G(n-1,k-1)+{n-k\choose k}t_{n-2k}.\]
Hence we have
\[G(n,k)=\sum_{j=0}^{k}{n-k\choose k-j}t_{n-2k+j}=\sum_{j=0}^k{n-k\choose 
j}t_{n-k-j}.\]

\smallskip

But if we have an explicit formula for $G$ then we have one for $f$ too, in 
view of (\ref{eq:gfrel}).  Indeed if we subtract (\ref{eq:gfrel}) with $k$ 
replaced by $k+1$ from (\ref{eq:gfrel}) we find that
\begin{eqnarray*}
f(n,k+1)&=&G(n,k)-G(n,k+1)\\
&=&\sum_{j=0}^k{n-k\choose k-j}t_{n-2k+j}-\sum_{j=0}^{k+1}{n-k-1\choose 
k+1-j}t_{n-2k-2+j}
\, ,
\end{eqnarray*}
thus proving our claim. \bx

\smallskip

Now if we use the asymptotic formula (\ref{eq:invol}) 
it is easy to see, from (\ref{eq:fform}), that
\[\lim_{n\to\infty}\frac{f(n,k)}{t_n}=\frac{k-1}{k!}\, ,\]
which we had previously derived from Corollary \ref{cor:first}. But now we can
use the asymptotic formula with a little more detail, on (\ref{eq:fform}), and
obtain the rate of approach to the limit.
\begin{theorem}
For each $k=3,4\dots$, the probability that
$k$ occurs in the $(1,2)$ position of a Young tableau of $n$ cells 
is 
\begin{equation}
\label{eq:secondterm}
     \frac{k-1}{k!} + \frac{k-4}{3(k-3)!}\, \frac{1}{n^{3/2}} + O\left(\frac{1}{n^2}\right) \qquad (n\to\infty).
\end{equation}
\end{theorem}

\end{document}